\documentclass[12pt]{article}
\usepackage{a4}
\usepackage{amsfonts}
\usepackage{amssymb}

\begin{document}
\title{\textbf{Determination of the real poles of the Igusa zeta function for curves}}
\author{Denis Ibadula \and Dirk Segers\thanks{Postdoctoral Fellow of the
Fund for Scientific Research - Flanders (Belgium).
}}

\date{June 23, 2010}
\maketitle

\begin{abstract}
The numerical data of an embedded resolution determine the
candidate poles of Igusa's $p$-adic zeta function. We determine in
complete generality which real candidate poles are actual poles in
the curve case.
\end{abstract}

\section{Introduction}

Several mathematicians have already obtained partial results about
the determination of the poles of Igusa's $p$-adic zeta function
for curves. In this paper, we will determine the real poles for an
arbitrary polynomial $f$ in two variables which is defined over a
$p$-adic field. People are interested in the poles of Igusa's
$p$-adic zeta function $Z_f(s)$ because they determine the
asymptotic behaviour of the number of solutions of polynomial
congruences and because they are the subject of the monodromy
conjecture (see for example \cite{Denefreport}).

Historically, one considered first only curves which are
absolutely analytically irreducible. Partial results were obtained
by Igusa \cite{Igusafirstterm} and Strauss \cite{Strauss}. Meuser
\cite{Meuser} determined the real poles, but she did not consider
the candidate pole $-1$. In 1985 Igusa \cite{Igusacomplexpowers}
solved that problem completely. He proved that the candidate poles
associated to the strict transform of $f$ are poles when the
domain of integration is small enough. Moreover, another candidate
pole of the minimal embedded resolution of $f$ is a pole if and
only if it is associated to an exceptional curve which is
intersected by three other irreducible components of the pull-back
of $f$. We have incorporated a generalization of this result
(Proposition 2).

In the general case, Loeser \cite{Loeser} obtained that an
exceptional curve $E_i$ does not contribute to the poles of
$Z_f(s)$ if $E_i$ is intersected one or two times by other
components of the pull-back of $f$ and if there are no other
intersection points over an algebraic closure. This was first
proved by Strauss in the absolutely analytically irreducible case,
where the last condition is automatically satisfied.

The next paper we want to mention is \cite{VeysIgusaCurves} of
Veys. He considers a polynomial $f$ in two variables over a number
field $F$ and takes the minimal embedded resolution of $f$ over an
algebraic closure of $F$. This setup allowed him to use a formula
\cite{Denefdegree} of Denef for $Z_f(s)$, which is valid for
almost all $p$-adic completions of $F$. He supposes that all
intersection points of irreducible components of the pull-back of
$f$ are defined over $F$. Under this condition, he proves the
converse of the result of Loeser for real candidate poles and for
almost all $p$-adic completions of $F$. Moreover, he deals with
the problem of a possible cancellation of several contributions to
the same real candidate pole.

In the proofs of the mentioned vanishing and non-vanishing
results, one needed certain relations between the various
numerical data of the embedded resolution. They were
systematically derived in \cite{Strauss}, \cite{Meuser} and
\cite{Igusacomplexpowers} for absolutely analytically irreducible
curves and finally, Loeser \cite{Loeser} obtained the necessary
relations in the general case. Igusa \cite{Igusacomplexpowers} and
Loeser \cite{Loeser} used a formula of Langlands \cite{Langlands}
to calculate the contribution of an exceptional curve to the
residue of $Z_f(s)$ at a candidate pole of candidate order one. We
will use a slight variant of this formula which was obtained in
\cite{Segersmathz}. Given an embedded resolution written as a
composition of blowing-ups, the second author explained there how
to calculate this contribution to the residue at the stage where
the exceptional curve is created. In Proposition 1, we determine
when this contribution is zero and when not. For this, we need new
ideas. It is not at all a straightforward generalization of what
was already known. Finally in Section 4, we will prove that
contributions to the same candidate pole will not cancel out. For
this, we use that the dual embedded resolution graph is an ordered
tree. This was obtained in \cite{Veysdetermination} when the base
field is algebraically closed.

\vspace{0,5cm}

\noindent \textbf{Acknowledgements.} The first author is grateful
to Willem Veys for inviting her to an academic visit at the
University of Leuven in 2009.  During this visit the work to this
paper was started. We also want to thank Bart Bories and Willem
Veys for the conversations that we had with them concerning this
paper and for their remarks.

\section{Definitions and our tools}

Let $K$ be a $p$-adic field, i.e., an extension of $\mathbb{Q}_p$
of finite degree. Let $R$ be the valuation ring of $K$, $P$ the
maximal ideal of $R$ and $q$ the cardinality of the residue field
$R/P$. For $z \in K$, let $\mathrm{ord} \, z \in \mathbb{Z} \cup
\{+\infty\}$ denote the valuation of $z$ and $|z|=q^{-\mathrm{ord}
\, z}$ the absolute value of $z$.

Let $f(x_1,x_2) \in K[x_1,x_2]$ be a polynomial in two variables
over $K$ and put $x=(x_1,x_2)$. Let $X$ be an open and compact
subset of $K^2$. Igusa's $p$-adic zeta function of $f$ is defined
by
\[ Z_f(s)= \int_X |f(x)|^s \, |dx| \]
for $s \in \mathbb{C}$, $\mbox{Re}(s) > 0$, where $|dx|$ denotes
the Haar measure on $K^2$, so normalised that $R^2$ has measure
$1$. Igusa proved that it is a rational function of $q^{-s}$ by
calculating the integral on an embedded resolution of $f$.
Therefore, it extends to a meromorphic function $Z_f(s)$ on
$\mathbb{C}$ which is also called Igusa's $p$-adic zeta function
of $f$.

Let $g : Y \rightarrow X$ be an embedded resolution of $f$. Here,
$Y$ is a $K$-analytic manifold. The meaning of embedded resolution
in our context is explained in \cite[Section 3.2]{Igusabook}.
Write $g = g_1 \circ \cdots \circ g_t : Y = Y_t \rightarrow X=Y_0$
as a composition of blowing-ups $g_i : Y_i \rightarrow Y_{i-1}$,
$i \in T_e := \{1,\ldots,t\}$. The exceptional curve of $g_i$ and
also the strict transforms of this curve are denoted by $E_i$. The
closed submanifolds of $Y$ of codimension one which are the zero
locus of the strict transform of an irreducible factor of $f$ in
$K[x,y]$ are denoted by $E_j$, $j \in T_s$. The corresponding
transforms in $Y_i$, $i \in \{0,\ldots,t-1\}$, are denoted in the
same way. Note that we had to be careful with the notion of
irreducible, because $X$ is totally disconnected as a topological
space. Put $T=T_e \cup T_s$. For $i \in T$, let $N_i$ and
$\nu_i-1$ be the multiplicities of respectively $f \circ g$ and
$g^*dx$ along $E_i$. The $(N_i,\nu_i)$ are called the numerical
data of $E_i$.

Let us recall Igusa's proof of the rationality of $Z_f(s)$. As we
already said, we calculate the defining integral on $Y$:
\[
Z_f(s) = \int_X |f(x)|^s \, |dx| = \int_Y |f \circ g|^s \,
|g^*dx|.
\]
Let $b$ be an arbitrary point of $Y$. There are three cases. In
the first case, there are two varieties $E_i$ and $E_j$, with $i,j
\in T$, that pass through $b$. We take a neighborhood $V$ of $b$
and analytic coordinates $(y_1,y_2)$ on $V$ such that $y_1$ is an
equation of $E_i$, $y_2$ is an equation of $E_j$,
\[
f \circ g = \varepsilon y_1^{N_i} y_2^{N_j} \qquad \mbox{and} \qquad
g^* dx = \eta y_1^{\nu_i-1} y_2^{\nu_j-1} dy
\]
on $V$ for non-vanishing $K$-analytic functions $\varepsilon$ and
$\eta$ on $V$. We may suppose that $y(V)=P^{k_1} \times P^{k_2}$,
with $k_1,k_2 \in \mathbb{Z}_{\geq 0}$, and that $|\varepsilon|$
and $|\eta|$ are constant on $V$. We get
\begin{eqnarray*}
\int_V |f \circ g|^s \, |g^* dx| & = & \int_{P^{k_1} \times
P^{k_2}} |\varepsilon|^s |\eta| |y_1|^{N_is+\nu_i-1}
|y_2|^{N_js+\nu_j-1} \, |dy| \\ & = & |\varepsilon|^s |\eta|
\left( \frac{q-1}{q} \right)^2
\frac{q^{-k_1(N_is+\nu_i)}}{1-q^{-(N_is+\nu_i)}}
\frac{q^{-k_2(N_js+\nu_j)}}{1-q^{-(N_js+\nu_j)}}.
\end{eqnarray*}
Note that this is a rational function of $q^{-s}$. In the second
case, there is one variety $E_i$, $i \in T$, that passes through
$b$. We take a neighborhood $V$ of $b$ and analytic coordinates
$(y_1,y_2)$ on $V$ such that $y_1$ is an equation of $E_i$,
\[
f \circ g = \varepsilon y_1^{N_i} \qquad \mbox{and} \qquad g^* dx
= \eta y_1^{\nu_i-1} dy
\]
on $V$ for non-vanishing $K$-analytic functions $\varepsilon$ and
$\eta$ on $V$. We may suppose that $y(V)=P^{k_1} \times P^{k_2}$,
with $k_1,k_2 \in \mathbb{Z}_{\geq 0}$, and that $|\varepsilon|$
and $|\eta|$ are constant on $V$. We get
\begin{eqnarray*}
\int_V |f \circ g|^s \, |g^* dx| & = & \int_{P^{k_1} \times
P^{k_2}} |\varepsilon|^s |\eta| |y_1|^{N_is+\nu_i-1} \, |dy| \\ &
= & |\varepsilon|^s |\eta| q^{-k_2} \frac{q-1}{q}
\frac{q^{-k_1(N_is+\nu_i)}}{1-q^{-(N_is+\nu_i)}}.
\end{eqnarray*}
In the third case, there is no variety $E_i$, $i \in T$, that
passes through $b$. We take a neighborhood $V$ of $b$ and analytic
coordinates $(y_1,y_2)$ on $V$ such that $f \circ g = \varepsilon$
and $g^* dx = \eta dy$ on $V$ for non-vanishing $K$-analytic
functions $\varepsilon$ and $\eta$ on $V$. We may suppose that
$y(V)=P^{k_1} \times P^{k_2}$, with $k_1,k_2 \in \mathbb{Z}_{\geq
0}$, and that $|\varepsilon|$ and $|\eta|$ are constant on $V$. We
get
\begin{eqnarray*}
\int_V |f \circ g|^s \, |g^* dx| & = & |\varepsilon|^s |\eta|
q^{-k_1-k_2}.
\end{eqnarray*}
It follows now that $Z_f(s)$ is a rational function of $q^{-s}$
because we can partition $Y$ into sets $V$ of the above form.

We obtain also from this calculation that every pole of $Z_f(s)$
is of the form
\[
-\frac{\nu_i}{N_i} + \frac{2k \pi \sqrt{-1}}{N_i \log q},
\]
with $k \in \mathbb{Z}$ and $i \in T$. These values are called the
candidate poles of $Z_f(s)$. If $i \in T$ is fixed, the values
$-\nu_i/N_i + (2k \pi \sqrt{-1})/(N_i \log q)$, $k \in
\mathbb{Z}$, are called the candidate poles of $Z_f(s)$ associated
to $E_i$. Because the poles of $1/(1-q^{-N_is- \nu_i})$ have order
one, we define the expected order of a candidate pole $s_0$ as the
highest number of $E_i$'s with candidate pole $s_0$ and with
non-empty intersection. The order of $s_0$ is of course less than
or equal to its expected order and a candidate pole $s_0$ of
expected order one is a pole if and only if the residue of
$Z_f(s)$ at $s_0$ is different from 0.

Let us explain the formula for the residue that we will use. Let
$s_0$ be a candidate pole of $E_i$, $i \in T$, and suppose that
$s_0$ is not a candidate pole of any $E_j$, with $j \in T$ and $j
\not= i$, which intersects $E_i$ in $Y$. Let $U$ be an open and
compact subset of $E_i$. The contribution of $U$ to the residue of
$Z_f(s)$ at $s_0$ is by definition the contribution to the residue
of $Z_f(s)$ at $s_0$ of an open and compact subset $V$ of $Y$ which
satisfies $V \cap E_i = U$ and which is disjoint from every other
$E_j$ with candidate pole $s_0$. Suppose that $U$ already exists in
$Y_r$ and if $i \in T_s$ we also suppose that it is non-singular in
$Y_r$. Suppose that $W$ is an open and compact subset of $Y_r$ for
which $W \cap E_i = U$ and that $(z_1,z_2)$ are analytic coordinates
on $W$ such that $z_1=0$ is an equation of $U$ on $W$. Write
\[
f \circ g_1 \circ \cdots \circ g_r = \gamma z_1^{N_i} \qquad
\mbox{and} \qquad (g_1 \circ \cdots \circ g_r)^* dx = \delta
z_1^{\nu_i-1} dy
\]
on $W$, for $K$-analytic functions $\gamma$ and $\delta$ on $W$.
Then, the contribution of $U$ to the residue of $Z_f(s)$ at $s_0$
is equal to
\begin{eqnarray} \label{residue}
\frac{q-1}{q N_i \log q} \left[ \int_{U} |\gamma|^s |\delta|
|dz_2| \right]^{\mathrm{mc}}_{s=s_0},
\end{eqnarray}
where $[ \cdot ]^{\mathrm{mc}}_{s=s_0}$ denotes the evaluation in
$s=s_0$ of the meromorphic continuation of the function between
the brackets. This formula was obtained by Langlands
\cite{Langlands} in the case $r=t$ and in general by the second
author in \cite{Segersmathz}.

We explain now the relations that we will need. Fix $r \in T_e$.
The exceptional curve $E_r$ is obtained by blowing-up at a point
$P \in Y_{r-1}$. Let $y=(y_1,y_2)$ be local coordinates on
$Y_{r-1}$ centered at $P$. Write in these local coordinates
\[
f \circ g_1 \circ \cdots \circ g_{r-1} = d \left( \prod_{i \in S}
(a_{i2}y_1-a_{i1}y_2)^{M_i} \right) \left( \prod_{i \in S'}
h_i^{M_i}(y_1,y_2) \right) + \mbox{ terms of higher degree},
\]
where all factors $a_{i2}y_1-a_{i1}y_2$ and $h_i$ are essentially
different (i.e. no factor is equal to another multiplied by an
element of $K^{\times}$) polynomials over $K$, where the $h_i$ are
irreducible homogeneous polynomials of degree at least two, where
$M_i \geq 1$ for every $i \in S \cup S'$ and where $d \in
K^{\times}$. Write also
\[
(g_1 \circ \cdots \circ g_{r-1})^*dx = \left( e \prod_{i \in S}
(a_{i2}y_1-a_{i1}y_2)^{\mu_i-1} + \mbox{terms of higher degree}
\right) dy,
\]
where $\mu_i \geq 1$ for every $i \in S$ and $e \in K^{\times}$.
Let $s_0=-\nu_r/N_r + (2k\pi \sqrt{-1})/(N_r \log q)$ be an
arbitrary candidate pole of $Z_f(s)$ associated to $E_r$. We
advise the reader to specialize everything what follows in this
section to the case $k=0$. Put $\alpha_i:=\mu_i+s_0M_i$ for every
$i \in S$. Because
\[
N_r = \sum_{i \in S} M_i + \sum_{i \in S'} (\deg h_i) M_i  \qquad
\mbox{and} \qquad \nu_r = \sum_{i \in S} (\mu_i - 1) + 2,
\]
it is straightforward to check that
\begin{eqnarray} \label{relation}
\sum_{i \in S} (\alpha_i -1) + \sum_{i \in S'} s_0(\deg h_i) M_i =
-2 + \frac{2k\pi \sqrt{-1}}{\log q}.
\end{eqnarray}
We now give another description of the $\alpha_i$. Let $F_i$ be
the point on $E_r$ which has coordinates $(a_{i1}:a_{i2})$ with
respect to the homogenous coordinates $(y_1:y_2)$ on $E_r \subset
Y_r$. Let $j$ be the unique element of $T \setminus \{r\}$ such
that $E_j$ passes through $F_i$ in $Y$. Let $\rho$ be the number
of blowing-ups among $g_r,\ldots,g_t$ which are centered at $F_i$.
Then, the announced description is $\alpha_i = \nu_j + s_0 N_j -(2
\rho k\pi \sqrt{-1})/(\log q)$. The second author proved this in
\cite[Section 2.7]{Segersmathz} in the case $k=0$, and the general
case is treated in a similar way. It follows that
$\mbox{Re}(\alpha_i)<0$ if and only if $-\nu_r/N_r < -\nu_j/N_j$.
One checks also easily that
\begin{eqnarray*}
s_0 \mbox{ is a candidate pole of } E_j & \Longleftrightarrow &
\nu_j + s_0 N_j \mbox{ is a multiple of } 2\pi \sqrt{-1}/(\log q)
\\ & \Longleftrightarrow & \alpha_i \mbox{ is a multiple of } 2
\pi \sqrt{-1}/(\log q).
\end{eqnarray*}
It is proved in \cite[Proposition II.3.1]{Loeser} that
$\mbox{Re}(\alpha_i) < 1$. Together with $(\ref{relation})$, this
implies that $\mbox{Re}(\alpha_i) \geq -1$ and that there is at
most one $i \in S$ with $\mbox{Re}(\alpha_i) < 0$.

\vspace{-0,2cm}

\section{Contribution of one exceptional curve}

\noindent \textsl{Setting of Proposition 1 and 2.} Let $f \in
K[x_1,x_2]$ and let $X$ be an open and compact subset of $K^2$.
Let $g:Y \rightarrow X$ be an embedded resolution of $f$. Write
$g=g_1 \circ \cdots \circ g_t: Y=Y_t \rightarrow X=Y_0$ as a
composition of blowing-ups $g_i: Y_i \rightarrow Y_{i-1}$, $i \in
T_e := \{1,\ldots,t\}$. The exceptional curve of $g_i$ and also
the strict transforms of this curve are denoted by $E_i$. Let $r
\in T_e$. The exceptional curve $E_r$ is obtained by blowing-up at
a point $P \in Y_{r-1}$. Let $(y_1,y_2)$ be local coordinates on
$Y_{r-1}$ centered at $P$. Write in these local coordinates
\[
f \circ g_1 \circ \cdots \circ g_{r-1} = d \left( \prod_{i \in S}
(a_{i2}y_1-a_{i1}y_2)^{M_i} \right) \left( \prod_{i \in S'}
h_i^{M_ i}(y_1,y_2) \right) + \mbox{ terms of higher degree},
\]
where all factors $a_{i2}y_1-a_{i1}y_2$ and $h_i$ are essentially
different polynomials over $K$, where the $h_i$ are irreducible
homogeneous polynomials of degree at least two, where $M_i \geq 1$
for every $i \in S \cup S'$ and where $d \in K^{\times}$. Write also
\[
(g_1 \circ \cdots \circ g_{r-1})^*dx = \left( e \prod_{i \in S}
(a_{i2}y_1-a_{i1}y_2)^{\mu_i-1} + \mbox{terms of higher degree}
\right) dy,
\]
where $\mu_i \geq 1$ for every $i \in S$ and $e \in K^{\times}$.

\vspace{0,4cm}

\noindent \textbf{Proposition 1.} Let $s_0:=-\nu_r/N_r$ be the
real candidate pole of $Z_f(s)$ associated to $E_r$. Suppose that
$\alpha_i:=\mu_i+s_0M_i \not= 0$ for every $i \in S$. Let
$\mathcal{R}$ be the contribution of $E_r$ to the residue of
$Z_f(s)$ at $s_0$. Then, $\mathcal{R} \not= 0$ if and only if $|S|
\geq 3$ or $|S'| \geq 1$. Moreover, if $\mathcal{R} \not= 0$, then
\begin{enumerate} \vspace{-0,1cm}
\item $\mathcal{R} > 0$ if and only if $\alpha_i > 0$ for every $i
\in S$ and \vspace{-0,1cm}
\item $\mathcal{R} < 0$ if and only if $\alpha_i < 0$ for some (and
thus exactly one) $i \in S$.
\end{enumerate}

\vspace{0,5cm}

\noindent \textbf{Proof.} If the number of elements of $S$ is one
or two and $S'$ is empty, then it is well know that
$\mathcal{R}=0$. We mentioned already in the introduction that
Loeser \cite{Loeser} proved this by using the formula of
Langlands, and the second author proved this again \cite[Section
3.1]{Segersmathz} as an illustration of his variant of this
formula.

Suppose from now on that $|S| \geq 3$ or $|S'| \geq 1$. We
consider first the case in which there exists one element $l \in
S$ satisfying $\alpha_l < 0$. Denote $Q:=S \setminus \{l\}$ and
$Q'=S'$. By applying an affine coordinate transformation, we may
assume that
\[
f \circ g_1 \circ \cdots \circ g_{r-1} = d \left( y_2^{M_l}
\prod_{i \in Q} (y_1-a_iy_2)^{M_i} \right) \left( \prod_{i \in Q'}
h_i^{M_ i}(y_1,y_2) \right) + \mbox{ terms of higher degree},
\]
and
\[
(g_1 \circ \cdots \circ g_{r-1})^*dx = \left( ey_2^{\mu_l-1}
\prod_{i \in Q} (y_1-a_iy_2)^{\mu_i-1} + \mbox{terms of higher
degree} \right) dy_,
\]
where the $a_i$, $i \in Q$, are different elements of $R$, where
$h_i$, $i \in Q'$, are different irreducible homogeneous
polynomials over $R$ of degree $d_i \geq 2$ with coefficient of
$y_1^{d_i}$ equal to 1, where $M_i \geq 1$ for every $i \in Q \cup
Q'$ and where $d,e \in K^{\times}$. We have that $\mathcal{R}$ is
the sum of two contributions, which we calculate on two different
charts by using formula (\ref{residue}). For the first
contribution, we look at the coordinates $(z_1,z_2)$ on $Y_r$ for
which $g_r(z_1,z_2) = (z_1,z_1z_2)$, and obtain $\kappa :=
(q-1)/(qN_r \log q)$ times
\begin{eqnarray*}
\lefteqn{\left[ |d|^s |e| \int_{P} |z_2|^{M_ls+\mu_l-1} \prod_{i
\in Q} \left( |1-a_iz_2|^{M_is+\mu_i-1} \right) \left( \prod_{i
\in Q'}
|h_i(1,z_2)|^{M_is} \right) \, |dz_2| \right]^{mc}_{s=s_0}} \hspace{4cm} \\
& = & |d|^{s_0} |e| \left[ \int_{P} |z_2|^{M_ls+\mu_l-1} \, |dz_2|
\right]^{mc}_{s=s_0} \hspace{4cm}
\\ & = & |d|^{s_0} |e| \frac{q-1}{q} \frac{1}{q^{\alpha_l}-1}.
\end{eqnarray*}
For the second contribution, we look at the coordinates
$(z_1',z_2')$ on $Y_r$ for which $g_r(z_1',z_2') =
(z_1'z_2',z_2')$, and obtain $\kappa$ times
\[
\left[ |d|^s |e| \int_{R} \prod_{i \in Q} \left(
|z_1'-a_i|^{M_is+\mu_i-1} \right) \left( \prod_{i \in Q'}
|h_i(z_1',1)|^{M_is} \right) \, |dz_1'| \right]^{mc}_{s=s_0},
\]
which is according to Lemma 3 less than $\kappa$ times
\[
|d|^{s_0} |e| \left[ \int_{R} |z_1'-a|^{Ms+\mu-1} \, |dz_1'|
\right]^{mc}_{s=s_0} = |d|^{s_0} |e| \frac{q-1}{q}
\frac{1}{1-q^{-(Ms_0+\mu)}},
\]
where $a \in R$, $M:=\left( \sum_{i \in Q} M_i \right) +
\left(\sum_{i \in Q'} d_iM_i \right)$ and $\mu:= \left( \sum_{i
\in Q} (\mu_i-1) \right)+1$. By using that $\alpha_l +
(Ms_0+\mu)=0$, we obtain
\[
\frac{1}{q^{\alpha_l}-1} + \frac{1}{1-q^{-(Ms_0+\mu)}} = 0,
\]
and this implies that $\mathcal{R} < 0$.

Consider now the case in which $\alpha_i > 0$ for every $i \in S$.
This case is much easier. After calculating $\mathcal{R}$
analogously as in the previous case by using formula
(\ref{residue}), you see that $\mathcal{R}$ is a sum of positive
numbers and thus positive. You have to use that $|h|$ is a locally
constant function for an irreducible polynomial $h$ over $K$ in
one variable of degree at least 2. $\qquad \Box$

\vspace{0,5cm}

We still have to prove Lemma 3. First, we prove Lemma 2, which is
a special case of Lemma 3. In the proof of Lemma 2, we need Lemma
1.

\vspace{0,5cm}

\noindent \textbf{Lemma 1.} Let $h \in R[x]$ be an irreducible
monic polynomial of degree $d \geq 2$ in one variable. Then, there
exists a unique $r \in \mathbb{N}$ and an element $b$ in $R$ such
that
\[
|h(x)| = |(x-b)^d|  \qquad \mbox{ if } x \in R \mbox{ and } x
\not\equiv b \mbox{ mod } P^r
\]
and
\[
q^{-dr} \leq |h(x)| < q^{-d(r-1)} \qquad \mbox{ if } x \in R\mbox{
and } x \equiv b \mbox{ mod } P^r.
\]
Moreover, $b$ is only determined modulo $P^r$ and $|h(x)|$ is
constant on $b+P^r$.

\vspace{0,5cm}

\noindent \textbf{Proof.} Let $\beta_1,\ldots,\beta_d$ be the
roots of $h$ in an algebraic closure of $K$. Let
$L:=K(\beta_1,\ldots,\beta_d)$, let $R_L$ be the valuation ring of
$L$ and denote the extension of the norm on $K$ to $L$ also by $|
\cdot |$. Note that $\beta_1,\ldots,\beta_d$ are different because
we work in characteristic zero and that they are in $R_L$ because
$h$ is monic and $R_L$ is the integral closure of $R$ in $L$.
Because
\begin{eqnarray*}
|h(x)| & = & |(x-\beta_1)(x-\beta_2) \cdots (x-\beta_d)| \\
& = & |x-\beta_1| \, |x-\beta_2| \cdots |x-\beta_d|,
\end{eqnarray*}
we look at $|x-\beta_i|$.

Take $i \in \{1,\ldots,d\}$. Let $r_i$ be the largest natural
number for which there exists an element $b_i \in R$ satisfying
$|b_i-\beta_i| < q^{-(r_i-1)}$. Note that this largest natural
number exists because $\beta_i \not\in R$. Note also that $b_i$ is
only determined modulo $P^{r_i}$. One checks now that
\[
|x-\beta_i| = |x-b_i|  \qquad \mbox{ if } x \in R \mbox{ and } x
\not\equiv b_i \mbox{ mod } P^{r_i}
\]
and
\[
q^{-r_i} \leq |x-\beta_i| < q^{-(r_i-1)} \qquad \mbox{ if } x \in
R \mbox{ and } x \equiv b_i \mbox{ mod } P^{r_i}.
\]
Moreover, $|x-\beta_i|$ is constant on $b_i+P^{r_i}$.

Finally, we show that $r_1=r_2=\cdots=r_d$ and we will denote this
by $r$. We also show that the $b_i$, $i \in \{1,\ldots,d\}$, are
the same modulo $P^r$, and because the $b_i$ are only determined
modulo $P^r$, we can take $b_1=b_2=\cdots=b_d$ and we will denote
this by $b$. Then the lemma obviously follows.

Take $i,j \in \{1,\ldots,d\}$ and suppose that $r_i \geq r_j$.
Because $\beta_i$ and $\beta_j$ are conjugate and because $b_i \in
R$, we get that also $b_i-\beta_i$ and $b_i-\beta_j$ are
conjugate. Because conjugate numbers have the same norm, we get
$|b_i-\beta_j| = |b_i-\beta_i| < q^{-(r_i-1)}$. This implies that
$r_i \leq r_j$ and thus $r_i=r_j$ and that $b_i \equiv b_j \mbox{
mod } P^{r_i}$. $\qquad \Box$

\vspace{0,5cm}

\noindent \textbf{Lemma 2.} Let $s_0$ be a negative rational
number. Let $a_1,\ldots,a_k$ be different elements of $R$ and
suppose that we have for every $i \in \{1,\ldots,k\}$ integers
$M_i \geq 1$ and $\mu_i \geq 1$ satisfying $\alpha_i :=
\mu_i+s_0M_i < 1$. Let $h_{k+1},\ldots,h_l \in R[x]$ be different
irreducible monic polynomials in one variable of degree at least
two. Denote the degree of $h_i$ by $d_i$ and suppose that we have
for every $i \in \{ k+1,\ldots,l\}$ an integer $M_i \geq 1$. Let
$r_i$ be the natural number which we associated to $h_i$ in the
previous lemma, and $b_i$ a corresponding element of $R$ which is
determined modulo $P^{r_i}$. Suppose that $r_{k+1}= \cdots = r_l$
and denote this number by $r$. Suppose also that $a_1 \equiv
\cdots \equiv a_k \equiv b_{k+1} \equiv \cdots \equiv b_l \mbox{
mod } P^r$ and that $a_i \not\equiv a_j \mbox{ mod } P^{r+1}$ for
$i \not= j$. Take now $a \in R$ such that $a \equiv a_1 \mbox{ mod
} P^r$ and put $M:=M_1+\cdots+M_k+d_{k+1}M_{k+1}+\cdots+d_lM_l$,
$\mu:= (\mu_1-1)+\cdots+(\mu_k-1)+1$ and $\alpha := \mu +s_0M$.
Suppose that $0<\alpha$ and that $k \geq 2$ or $l \geq k+1$. Then
\begin{eqnarray*}
\lefteqn{\left[ \int_{a+P^r} |x-a_1|^{M_1s+\mu_1-1} \cdots
|x-a_k|^{M_ks+\mu_k-1} |h_{k+1}(x)|^{M_{k+1}s} \ldots
|h_l(x)|^{M_ls} \, |dx| \right]^{mc}_{s=s_0}} \hspace{4cm} \\ & <
& \left[ \int_{a+P^r} |x-a|^{Ms+\mu-1} \, |dx|
\right]^{mc}_{s=s_0}. \hspace{4cm}
\end{eqnarray*}
Moreover, the integrands are the same for every $x \in R \setminus
(a+P^r)$.

\vspace{0,5cm}

\noindent \textsl{Remark.} (1) The conditions $a_i \equiv a_j
\mbox{ mod } P^r$ and $a_i \not\equiv a_j \mbox{ mod } P^{r+1}$
for $i,j \in \{1,\ldots,k\}$ with $i \not= j$ imply that $k \leq
q$. \\ (2) We have that
\[
\alpha - 1 = \sum_{i=1}^k (\alpha_i -1) + \sum_{i=k+1}^l s_0 d_i
M_i.
\]
Consequently, the condition $\alpha_i < 1$ for every $i \in
\{1,\ldots,k\}$ implies that $\alpha < \alpha_i$.

\vspace{0,5cm}

\noindent \textbf{Proof.} In the first step, we reduce to the case
in which the polynomials $h_i$ do not occur. We have
\begin{eqnarray*}
\lefteqn{\left[ \int_{a+P^r} |x-a|^{Ms+\mu-1} \, |dx|
\right]^{mc}_{s=s_0}} \hspace{1cm}\\ & = & \frac{q-1}{q}
\frac{q^{-r\alpha}}{1-q^{-\alpha}} \\ & = & \frac{q-1}{q}
\frac{q^{-r((M_1+\cdots+M_k)s_0+\mu)}}{1-q^{-\alpha}}
q^{-r(d_{k+1}M_{k+1}+\cdots+d_lM_l)s_0} \\ & \geq & \frac{q-1}{q}
\frac{q^{-r((M_1+\cdots+M_k)s_0+\mu)}}{1-q^{-((M_1+\cdots+M_k)s_0+\mu)}}
q^{-r(d_{k+1}M_{k+1}+\cdots+d_lM_l)s_0} \\ & = &
q^{-r(d_{k+1}M_{k+1}+\cdots+d_lM_l)s_0} \left[ \int_{a+P^r}
|x-a|^{(M_1+\cdots+M_k)s+\mu-1} \, |dx| \right]^{mc}_{s=s_0},
\end{eqnarray*}
with a strict inequality if $l \geq k+1$. Because $|h_i|$ is
constant on $a+P^r$ with $q^{-d_ir} \leq |h_i(a)|$, $s_0<0$ and
the second factor on the right hand side of the following
inequality is positive (this follows from the calculation of this
factor in the second part of the proof), we get
\begin{eqnarray*}
\lefteqn{\left[ \int_{a+P^r} |x-a_1|^{M_1s+\mu_1-1} \cdots
|x-a_k|^{M_ks+\mu_k-1} |h_{k+1}(x)|^{M_{k+1}s} \ldots
|h_l(x)|^{M_ls} \, |dx| \right]^{mc}_{s=s_0}} \\ & \leq &
q^{-r(d_{k+1}M_{k+1}+\cdots+d_lM_l)s_0} \left[ \int_{a+P^r}
|x-a_1|^{M_1s+\mu_1-1} \cdots |x-a_k|^{M_ks+\mu_k-1} \, |dx|
\right]^{mc}_{s=s_0}.
\end{eqnarray*}
These two inequalities imply that it is enough to consider the
case where the $h_i$ do not occur.

So in the second step, we prove that
\[
\left[ \int_{a+P^r} |x-a|^{Ms+\mu-1} \, |dx| \right]^{mc}_{s=s_0}
> \left[ \int_{a+P^r} |x-a_1|^{M_1s+\mu_1-1} \cdots
|x-a_k|^{M_ks+\mu_k-1} \, |dx| \right]^{mc}_{s=s_0}
\]
if $k \geq 2$, where $M=M_1+\cdots+M_k$. We calculate both sides.
We partition the domain of integration of the integral on the
right hand side into the following $k+1$ sets: $a_1+P^{r+1},
\ldots ,a_k+P^{r+1}$ and the set consisting of all other points of
$a+P^r$. In this way, the above inequality becomes
\begin{eqnarray*}
\frac{q-1}{q} \frac{q^{-(r-1)\alpha}}{q^{\alpha}-1} & > &
\frac{q-1}{q} \frac{q^{-r\alpha_1}}{q^{\alpha_1}-1}
q^{-r(\alpha_2-1)-r(\alpha_3-1)-\ldots-r(\alpha_k-1)} + \cdots \\
&  &  + \frac{q-1}{q} \frac{q^{-r\alpha_k}}{q^{\alpha_k}-1}
q^{-r(\alpha_1-1)-r(\alpha_2-1)-\ldots-r(\alpha_{k-1}-1)} \\ &  &
+ \frac{q-k}{q^{r+1}}
q^{-r(\alpha_1-1)-r(\alpha_2-1)-\ldots-r(\alpha_k-1)}.
\end{eqnarray*}
By using the fact that $\alpha-1 = \sum_{i=1}^k ( \alpha_i-1)$,
this is equivalent to
\[
(q-1) \frac{q^{\alpha}}{q^{\alpha}-1} > \frac{q-1}{q^{\alpha_1}-1}
+ \cdots + \frac{q-1}{q^{\alpha_k}-1} + q-k
\]
and thus also to
\[
\frac{1}{q^{\alpha}-1} + \frac{k-1}{q-1} >
\frac{1}{q^{\alpha_1}-1} + \cdots + \frac{1}{q^{\alpha_k}-1}.
\]
Consider the function
\[
h : \, ]0,1] \rightarrow \mathbb{R} : x \mapsto \frac{1}{q^x-1}.
\]
One checks easily that $h$ is convex, i.e. $h''(x) > 0$ for every
$x \in ]0,1[$. Consider the linear function $g$, i.e. polynomial
function of degree one, determined by $g(\alpha) = h(\alpha) =
1/(q^{\alpha}-1)$ and $g(1)=h(1)=1/(q-1)$. Then
\begin{eqnarray*}
\frac{1}{q^{\alpha}-1} + \frac{k-1}{q-1} & = & g(\alpha) +
(k-1)g(1) \\ & = & g(\alpha_1) + \cdots + g(\alpha_k) \\ & > &
h(\alpha_1) + \cdots + h(\alpha_k) \\ & = &
\frac{1}{q^{\alpha_1}-1} + \cdots + \frac{1}{q^{\alpha_k}-1},
\end{eqnarray*}
where we used in the second line that $g$ is linear and that
$\alpha + k-1 = \alpha_1+\cdots+\alpha_k$ and in the third line
that $g$ is linear and $h$ convex, that $g(\alpha)=h(\alpha)$ and
$g(1)=h(1)$ and that $0<\alpha<\alpha_i<1$ for every $i \in
\{1,\ldots,k\}$.

The final statement in the lemma is easy. $\qquad \Box$

\vspace{0,5cm}

\noindent \textbf{Lemma 3.} Let $s_0$ be a negative rational
number. Let $\gamma,\delta \in R[x]$ be monic polynomials in one
variable. Suppose that $\delta$ factors into linear polynomials
over $R$ and that all roots of $\delta$ are also roots of
$\gamma$. Write
\[
\gamma(x) = \left( \prod_{i \in Q} (x-a_i)^{M_i} \right) \left(
\prod_{i \in Q'} h_i^{M_ i}(x) \right)
\]
where the $a_i$, $i \in Q$, are different elements of $R$, where
$h_i$, $i \in Q'$, are different irreducible monic polynomials
over $R$ of degree at least two, and where $M_i \geq 1$ for every
$i \in Q \cup Q'$. Denote the degree of $h_i$ by $d_i$. Write also
\[
\delta(x) = \prod_{i \in Q} (x-a_i)^{\mu_i-1},
\]
where $\mu_i \geq 1$ for every $i \in Q$. Take any $a \in R$ and
put $M:=\left( \sum_{i \in Q} M_i \right) + \left(\sum_{i \in Q'}
d_iM_i \right)$ and $\mu:= \left( \sum_{i \in Q} (\mu_i-1)
\right)+1$. Suppose that $0 < \alpha := \mu + s_0 M$ and $1
> \alpha_i := \mu_i +s_0 M_i$ for every $i \in Q$ and that $|Q| \geq 2$ or $|Q'|
\geq 1$. Then
\[
\left[ \int_{R} |\gamma(x)|^s |\delta(x)| \, |dx|
\right]^{mc}_{s=s_0} < \left[ \int_{R} |x-a|^{Ms+\mu-1} \, |dx|
\right]^{mc}_{s=s_0}.
\]

\vspace{0,5cm}

\noindent \textbf{Proof.} We associate a tree to a monic
polynomial $g \in R[x]$ in one variable as follows. If $a_1, a_2
\in R$ are roots of $g$ and if $a_1 \equiv a_2 \mbox{ mod } P^r$
and $a_1 \not\equiv a_2 \mbox{ mod } P^{r+1}$, then we associate a
bullet to $a_1 + P^r$. If we have an irreducible factor of $g$ of
degree at least two, we have associated $r \in \mathbb{N}$ and $b
+ P^r$ to it, and we associate to $b + P^r$ a bullet. (If say $a
+P^r$ occurs several times in this way, we associate to it only
one bullet.) We make a tree by connecting bullets in the obvious
way, i.e. if we have bullets associated to $a + P^r$ and $b +
P^t$, with $r > t$, then we connect both bullets if $a+P^r \subset
b + P^t$ and for any $c+P^u$ corresponding to another bullet we
don't have $a+P^r \subset c+P^u \subset b + P^t$. Note that this
tree is finite and has one root.

We start at the left hand side of the inequality that we want to
prove and we consider the tree associated to $\gamma$. We will
construct step by step other integrands for which the associated
tree is the previous one except one bullet at the end of the tree.
So the tree becomes easier after each step. We will do this until
the tree completely disappears. At this stage, the integrand will
be the one on the right hand side of the inequality that we want
to prove.

We explain the first step. We take a bullet at the end of the tree
of $\gamma$. This bullet is associated to an element of $R/P^r$
for some $r$, let us say $a_0 + P^r$. Let $a_1,\ldots,a_k$ be all
the roots of $\gamma$ for which $a_i \equiv a_0 \mbox{ mod } P^r$.
Note that $a_i \not\equiv a_j \mbox { mod } P^{r+1}$ for $i,j \in
\{1,\ldots,k\}$ with $i \not= j$, because we took a bullet at the
end of the tree. Let $h_{k+1},\ldots,h_l$ be all the irreducible
factors of $\gamma$ to which we have associated $a_0 + P^r$. Write
\[
\gamma(x) = \widetilde{\gamma}(x) (x-a_1)^{M_1} \cdots
(x-a_k)^{M_k} h_{k+1}^{M_{k+1}} \cdots h_l^{M_l}
\]
and
\[
\delta(x) = \widetilde{\delta}(x) (x-a_1)^{\mu_1-1} \cdots
(x-a_k)^{\mu_k-1}.
\]
Put $M_0:=M_1+\cdots+M_k+d_{k+1}M_{k+1}+\cdots+d_lM_l$ and
$\mu_0:= (\mu_1-1)+\cdots+(\mu_k-1)+1$. Put $\gamma_1(x) =
\widetilde{\gamma}(x) (x-a_0)^{M_0}$ and $\delta_1(x) =
\widetilde{\delta}(x) (x-a_0)^{\mu_0-1}$. We have
\begin{eqnarray*}
\lefteqn{\left[ \int_{a_0+P^r} |\gamma(x)|^s |\delta(x)| \, |dx|
\right]^{mc}_{s=s_0}} \\ & = & \left[ |\widetilde{\gamma}(a_0)|^s
|\widetilde{\delta}(a_0)| \int_{a_0+P^r} \prod_{i=1}^k
|x-a_i|^{M_is+\mu_i-1} \prod_{j=k+1}^l |h_j(x)|^{M_js} \, |dx|
\right]^{mc}_{s=s_0} \\ & = & |\widetilde{\gamma}(a_0)|^{s_0}
|\widetilde{\delta}(a_0)| \left[ \int_{a_0+P^r} \prod_{i=1}^k
|x-a_i|^{M_is+\mu_i-1} \prod_{j=k+1}^l |h_j(x)|^{M_js} \, |dx|
\right]^{mc}_{s=s_0} \\ & < & |\widetilde{\gamma}(a_0)|^{s_0}
|\widetilde{\delta}(a_0)| \left[ \int_{a_0+P^r}
|x-a_0|^{M_0s+\mu_0-1} \, |dx| \right]^{mc}_{s=s_0}
\\ & = & \left[ |\widetilde{\gamma}(a_0)|^s
|\widetilde{\delta}(a_0)| \int_{a_0+P^r} |x-a_0|^{M_0s+\mu_0-1} \,
|dx| \right]^{mc}_{s=s_0} \\ & = & \left[ \int_{a_0+P^r}
|\widetilde{\gamma}(x)|^s |\widetilde{\delta}(x)|
|x-a_0|^{M_0s+\mu_0-1} \, |dx| \right]^{mc}_{s=s_0} \\ & = &
\left[ \int_{a_0+P^r} |\gamma_1(x)|^s |\delta_1(x)| \, |dx|
\right]^{mc}_{s=s_0},
\end{eqnarray*}
where we used the previous lemma and we used twice that
$|\widetilde{\gamma}|$ and $|\widetilde{\delta}|$ are constant on
$a_0+P^r$, and
\begin{eqnarray*}
\lefteqn{\left[ \int_{R \setminus (a_0+P^r)} |\gamma(x)|^s
|\delta(x)| \, |dx| \right]^{mc}_{s=s_0}} \\ & = & \left[ \int_{R
\setminus (a_0+P^r)} |\widetilde{\gamma}(x)|^s
|\widetilde{\delta}(x)| \prod_{i=1}^k |x-a_i|^{M_is+\mu_i-1}
\prod_{j=k+1}^l |h_j(x)|^{M_js} \, |dx| \right]^{mc}_{s=s_0} \\ &
= & \left[ \int_{R \setminus (a_0+P^r)} |\widetilde{\gamma}(x)|^s
|\widetilde{\delta}(x)| |x-a_0|^{M_0s+\mu_0-1} \, |dx|
\right]^{mc}_{s=s_0} \\ & = & \left[ \int_{R \setminus (a_0+P^r)}
|\gamma_1(x)|^s |\delta_1(x)| \, |dx| \right]^{mc}_{s=s_0},
\end{eqnarray*}
where we used the last sentence in the formulation of the previous
lemma, so that
\begin{eqnarray*}
\left[ \int_{R} |\gamma(x)|^s |\delta(x)| \, |dx|
\right]^{mc}_{s=s_0} & < & \left[ \int_{R} |\gamma_1(x)|^s
|\delta_1(x)| \, |dx| \right]^{mc}_{s=s_0}.
\end{eqnarray*}

In the second step, we do the same as in the first step, but use
now $\gamma_1(x)$ instead of $\gamma(x)$ and $\delta_1(x)$ instead
of $\delta(x)$. Remark that the $M$ and $\mu$ determined by
$\gamma$ and $\delta$ are the same as the analogous ones
determined by $\gamma_1$ and $\delta_1$. Remark also that the tree
associated to $\gamma_1$ is the tree associated to $\gamma$ with
one bullet missing.

Denote the number of bullets of the tree associated to $\gamma$ by
$w$. Then, after $w$ steps, the tree completely disappears. If the
root of the tree is associated to $a_0' + P^{r'}$, then
$\gamma_w(x)=(x-a_0')^M$ and $\delta_w(x)=(x-a_0')^{\mu-1}$.
Consequently,
\begin{eqnarray*}
\left[ \int_{R} |\gamma(x)|^s |\delta(x)| \, |dx|
\right]^{mc}_{s=s_0} & < & \left[ \int_{R} |\gamma_1(x)|^s
|\delta_1(x)| \, |dx| \right]^{mc}_{s=s_0} \\ & < & \left[
\int_{R} |\gamma_2(x)|^s |\delta_2(x)| \, |dx|
\right]^{mc}_{s=s_0} \\ & < & \ldots \\ & < & \left[ \int_{R}
|\gamma_w(x)|^s |\delta_w(x)| \, |dx| \right]^{mc}_{s=s_0} \\ & =
& \left[ \int_{R} |x-a_0'|^{Ms+\mu-1} \, |dx| \right]^{mc}_{s=s_0}
\\ & = & \left[ \int_{R} |x-a|^{Ms+\mu-1} \, |dx|
\right]^{mc}_{s=s_0}. \qquad \Box
\end{eqnarray*}

\vspace{0,5cm}

In the next proposition, we use the setting explained in the
beginning of this section.

\vspace{0,5cm}

\noindent \textbf{Proposition 2.} Let $s_0 := -\nu_r/N_r + (2k\pi
\sqrt{-1})/(N_r \log q)$ be an arbitrary candidate pole of
$Z_f(s)$ associated to $E_r$. Suppose that $\alpha_i :=
\mu_i+s_0M_i$ is not a multiple of $2\pi \sqrt{-1}/(\log q)$ for
every $i \in S$. Suppose $|S|=3$ and $|S'|=0$. Let $\mathcal{R}$
be the contribution of $E_r$ to the residue of $Z_f(s)$ at $s_0$.
Then, $\mathcal{R} \not= 0$.

\vspace{0,5cm}

\noindent \textbf{Proof.} Denote the elements of $S$ by $1$, $2$
and $3$. Then the equality (\ref{relation}) says $\alpha_1 +
\alpha_2 + \alpha_3 = 1 + (2k \pi \sqrt{-1})/(\log q)$. By
applying an affine coordinate transformation, we may assume that
\[
f \circ g_1 \circ \cdots \circ g_{r-1} = dy_2^{M_1} y_1^{M_2}
(y_1-ay_2)^{M_3} + \mbox{ terms of higher degree},
\]
and
\[
(g_1 \circ \cdots \circ g_{r-1})^* dx = \left( e y_2^{\mu_1-1}
y_1^{\mu_2-1} (y_1-ay_2)^{\mu_3-1} + \mbox{ terms of higher
degree} \right) dy
\]
with $a \in R \setminus P$ and $d,e \in K^{\times}$. We have now
that
\begin{eqnarray*}
\mathcal{R} & = & |d|^{s_0} |e| \left( \frac{q-1}{q}
\frac{1}{q^{\alpha_1}-1} + \frac{q-1}{q} \frac{1}{q^{\alpha_2}-1}
+ \frac{q-1}{q} \frac{1}{q^{\alpha_3}-1} + \frac{q-2}{q} \right) \\
& = & |d|^{s_0} |e| \left( \frac{1-q^{\alpha_1-1}}{1-q^{-\alpha_1}}
\cdot \frac{1-q^{\alpha_2-1}}{1-q^{-\alpha_2}} \cdot
\frac{1-q^{\alpha_3-1}}{1-q^{-\alpha_3}} \right) \\ & \not= & 0.
\end{eqnarray*}
The second equality can be checked by a straightforward
calculation and is due to Sally and Taibleson
\cite{SallyTaibleson}. $\qquad \Box$

\vspace{0,5cm}

\noindent \textsl{Remark.} (1) The determination of all poles
(real and complex) of an absolutely analytically irreducible curve
follows now immediately. This was one of the main results of the
paper \cite{Igusacomplexpowers} of Igusa. He also used
Sally-Taibleson's formula. \\ (2) However, if $|S| \not= 3$ or
$|S'| \not= 0$, it is not clear to us which non-real candidate
poles are poles and which are not. It can happen that a real
candidate pole is a pole, and that other candidate poles with the
same real part are not poles. This is for example the case when
$f=x_1^2+x_2^2$ and $p=2$ (see \cite[Example 2.8]{Segersmathz}).

\section{Poles of Igusa's $p$-adic zeta function}

Let $f \in K[x_1,x_2]$ and let $X$ be an open and compact subset of
$K^2$. Suppose that $f_{red}$ has only one singular point $P_0$ in
$X$. Let $g:Y \rightarrow X$ be an embedded resolution of $f$. Write
$g=g_1 \circ \cdots \circ g_t: Y=Y_t \rightarrow X=Y_0$ as a
composition of blowing-ups $g_i: Y_i \rightarrow Y_{i-1}$, $i \in
T_e := \{1,\ldots,t\}$, centered at $P_{i-1} \in Y_{i-1}$. The
exceptional curve of $g_i$ and also the strict transforms of this
curve are denoted by $E_i$. The closed submanifolds of $Y$ of
codimension one which are the zero locus of the strict transform of
an irreducible factor of $f$ in $K[x,y]$ are denoted by $E_j$, $j
\in T_s$. The corresponding transforms in $Y_i$, $i \in
\{0,\ldots,t-1\}$, are denoted in the same way. Put $T=T_e \cup
T_s$.

In the (dual) embedded resolution graph of the germ of $f$ at
$P_0$ one associates to each exceptional curve a vertex
(represented by a dot) and to each intersection between
exceptional curves in $Y$ an edge, connecting the corresponding
vertices. We also associate to each analytically irreducible
component of the strict transform of the germ of $f$ at $P_0$ a
vertex (represented by a circle), and to its (unique) intersection
with an exceptional curve in $Y$ a corresponding edge. It is clear
that this graph is a finite connected tree with all circles end
vertices.

Now to each vertex of the embedded resolution graph we associate
the corresponding ratio $\nu_i/N_i$. This makes the embedded
resolution graph into an ordered tree. More precisely, the
vertices for which the associated number is equal to $\min_{i \in
T} \nu_i/N_i$, together with their edges, form a connected part
$\mathcal{M}$ of the embedded resolution graph, and starting from
an end vertex of the minimal part $\mathcal{M}$, the numbers
$\nu_i/N_i$ strictly increase along any path in the tree (away
from $\mathcal{M}$). This follows from relation $(\ref{relation})$
and the bound on the $\alpha$'s, which imply for example that
there exists at most one $E_j$ which intersects a given $E_r$, $r
\in T_e$, in $Y$ with $\nu_j/N_j < \nu_r/N_r$ (see Section 2). For
more details, see \cite[Theorem 3.3]{Veysdetermination}, where the
base field is $\mathbb{C}$ instead of $K$, but nevertheless the
proof of our statement is similar.

\vspace{0,5cm}

\noindent \textsl{Example.} When $f$ is absolutely analytically
irreducible at $P_0$ with $g$ different Puiseux exponents, then
the resolution graph has the form

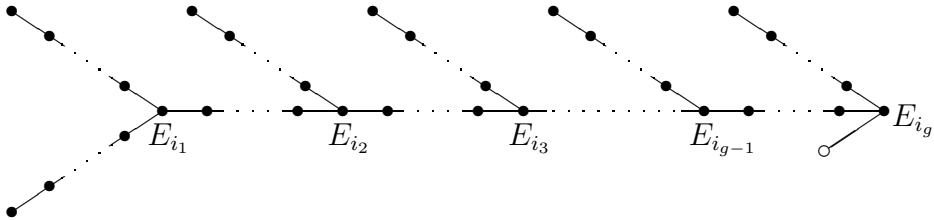
\begin{picture}(150,40)(0,0)
\put(0,33.333){\circle*{1,5}} \put(5,30){\circle*{1.5}}
\put(15,23.333){\circle*{1,5}}

\qbezier[6](5,30)(10,26.666)(15,23.333)
\put(0,33.333){\line(3,-2){7}} \put(20,20){\line(-3,2){7}}

\put(0,6.666){\circle*{1,5}} \put(5,10){\circle*{1,5}}
\put(15,16.666){\circle*{1,5}}

\qbezier[6](5,10)(10,13.333)(15,16.666)
\put(0,6.666){\line(3,2){7}} \put(20,20){\line(-3,-2){7}}

\put(20,20){\circle*{1,5}} \put(26,20){\circle*{1,5}}
\put(38,20){\circle*{1,5}}

\qbezier[6](26,20)(32,20)(38,20) \put(20,20){\line(1,0){8}}
\put(44,20){\line(-1,0){8}}

\put(24,33.333){\circle*{1,5}} \put(29,30){\circle*{1.5}}
\put(39,23.333){\circle*{1,5}}

\qbezier[6](29,30)(34,26.666)(39,23.333)
\put(24,33.333){\line(3,-2){7}} \put(44,20){\line(-3,2){7}}

\put(44,20){\circle*{1,5}} \put(50,20){\circle*{1,5}}
\put(62,20){\circle*{1,5}} \put(68,20){\circle*{1,5}}

\qbezier[6](50,20)(56,20)(62,20) \put(44,20){\line(1,0){8}}
\put(68,20){\line(-1,0){8}}

\put(48,33.333){\circle*{1,5}} \put(53,30){\circle*{1.5}}
\put(63,23.333){\circle*{1,5}}

\qbezier[6](53,30)(58,26.666)(63,23.333)
\put(48,33.333){\line(3,-2){7}} \put(68,20){\line(-3,2){7}}

\qbezier[12](68,20)(80,20)(92,20) \put(68,20){\line(1,0){3}}
\put(92,20){\line(-1,0){3}}

\put(72,33.333){\circle*{1,5}} \put(77,30){\circle*{1.5}}
\put(87,23.333){\circle*{1,5}}

\qbezier[6](77,30)(82,26.666)(87,23.333)
\put(72,33.333){\line(3,-2){7}} \put(92,20){\line(-3,2){7}}

\put(92,20){\circle*{1,5}} \put(98,20){\circle*{1,5}}
\put(110,20){\circle*{1,5}} \put(116,20){\circle*{1,5}}

\qbezier[6](98,20)(104,20)(110,20) \put(92,20){\line(1,0){8}}
\put(116,20){\line(-1,0){8}}

\put(96,33.333){\circle*{1,5}} \put(101,30){\circle*{1.5}}
\put(111,23.333){\circle*{1,5}} \put(108,14.666){\circle{1,5}}

\qbezier[6](101,30)(106,26.666)(111,23.333)
\put(96,33.333){\line(3,-2){7}} \put(116,20){\line(-3,2){7}}

\put(116,20){\line(-3,-2){7,375}}

\put(18,15.5){$E_{i_1}$} \put(42,15.5){$E_{i_2}$}
\put(66,15.5){$E_{i_3}$} \put(90,15.5){$E_{i_{g-1}}$}
\put(117,18){$E_{i_g}$}
\end{picture}

\noindent The minimal part $\mathcal{M}$ consists just of
$E_{i_1}$ (see \cite[Corollary 2.1]{Strauss} or \cite[Proposition
3.6]{Veysdetermination}).

\vspace{0,5cm}

\noindent \textbf{Theorem.} Suppose that we are in the setting of
the first paragraph of this section. Then,
\begin{enumerate}
\item a real number $s_0$ is a pole of order 2 if and only if
there exist $i,j \in T$ with $s_0=-\nu_i/N_i=-\nu_j/N_j$ such that
$E_i$ and $E_j$ intersect on $Y$, moreover, $Z_f(s)$ has at most
one real pole of order two, and if there is a pole of order two,
it is the pole closest to the origin,
\item a real number $s_0 \in \{ -\nu_i/N_i \mid i \in T_e \} \setminus
\{ -\nu_i/N_i \mid i \in T_s \}$ which is not a pole of order 2 is
a pole of order 1 if and only if there exists at least one $i \in
T_e$ with $s_0=-\nu_i/N_i$ such that $f \circ g_1 \circ \cdots
\circ g_{i-1}$ is given in local coordinates centered at $P_{i-1}$
by a power series with lowest degree part a homogeneous polynomial
which is not a power of a linear polynomial or a product of two
such powers, and
\item a real number $s_0 \in \{ -\nu_i/N_i \mid i \in T_s \}$ which
is not a pole of order 2 is a pole of order 1 for a small enough
open and compact neighborhood $X$ of $P_0$.
\end{enumerate}

\vspace{0,5cm}

\noindent \textbf{Proof.} (1) It is clear that there exist $i,j
\in T$ with $s_0=-\nu_i/N_i=-\nu_j/N_j$ such that $E_i$ and $E_j$
intersect on $Y$ if $s_0$ is a real pole of order 2. If $E_i$ and
$E_j$ intersect on $Y$ with $s_0=-\nu_i/N_i=-\nu_j/N_j$, then the
contribution of $P:=E_i \cap E_j$ to the coefficient $b_{-2}$ in
the Laurent series
\[
\frac{b_{-2}}{(s-s_0)^2} + \frac{b_{-1}}{s-s_0} + b_0 + b_1(s-s_0)
+ \cdots
\]
of $Z_f(s)$ at $s_0$ is equal to
\[
|\varepsilon(P)|^{s_0} |\eta(P)| \frac{(q-1)^2}{q^2 N_i N_j (\log
q)^2}
>0,
\]
and consequently, the contributions of different intersecting
pairs on $Y$ cannot cancel each other. The other statements follow
from the ordered tree structure of the embedded resolution graph.

(2) The `only if' part is the well known part of Proposition 1
which is due to Loeser. For the other implication, we have to use
Proposition 1 and the ordered tree structure of the embedded
resolution graph. There are two possibilities. In the first case,
$E_i$ is the minimal part of the embedded resolution graph. In
this case, there is only one contribution to the residue, which is
positive. In the other case, there is at least one non-zero
contribution to the residue, and all such contributions are
negative.

(3) There are two possibilities. In the first case, $s_0
=-\nu_i/N_i$, with $i \in T_s$ and $E_i$ is the minimal part of
the embedded resolution graph. In this case, there is only one
contribution to the residue, which is positive. In the other case,
we take a (small enough) open and compact neighborhood $V$ of
$\cup_{i \in T_e} E_i \subset Y$ on which all the $E_i$, $i \in
T_s$ with $s_0=-\nu_i/N_i$, have a negative contribution to the
residue of $Z_f(s)$ at $s_0$. The $E_i$, $i \in T_e$ with $s_0 =
-\nu_i/N_i$, have a non-positive contribution to the residue of
$Z_f(s)$ at $s_0$. If we replace $X$ by $g(V)$ or by an open and
compact neighborhood of $P$ contained in $g(V)$, we obtain what we
want. $\qquad \Box$

\vspace{0,5cm}

\noindent \textsl{Remark.} It follows from the previous theorem
and from the result of Loeser mentioned in the introduction that
$\mbox{Re}(s_0)$ is a pole of $Z_f(s)$ if $s_0$ is a pole of
$Z_f(s)$. Consequently, we know the set of real parts of poles of
$Z_f(s)$. This determines the asymptotic behaviour of the number
of solutions of the corresponding polynomial congruences (see
\cite{Segersasymptotic} for more information).

\footnotesize{

\vspace{0,5cm}

\noindent \textsc{Faculty of Mathematics and Informatics,
``Ovidius'' University, Constanta, Mamaia Bd. 124, RO 900527
Constanta, Romania}
\\ \textsl{E-mail address:} denis@univ-ovidius.ro

\vspace{0,5cm}

\noindent \textsc{University of Leuven, Department of
Mathematics, Celestijnenlaan 200B, B-3001 Leuven, Belgium} \\
\textsl{E-mail address:} dirk.segers@wis.kuleuven.be \\
\textsl{URL:} http://wis.kuleuven.be/algebra/segers/segers.htm
\end{document}